\documentclass[a4 paper, reqno]{amsart} 
\usepackage{titling}
\usepackage{amsaddr}
\usepackage{latexsym}
\usepackage{hyperref}
\hypersetup{
	colorlinks=true,
	linkcolor=red,
	citecolor=blue,
}
\usepackage{tikz-cd}
\usepackage{xcolor}
\usepackage
[
margin=3cm
]
{geometry}

\newtheorem{defn}{Definition}[section]
\newtheorem{thm}[defn]{Theorem}
\newtheorem{lemma}[defn]{Lemma}
\newtheorem{prop}[defn]{Proposition}
\newtheorem{corr}[defn]{Corollary}
\newtheorem{xmpl}[defn]{Example}
\newtheorem{rmk}[defn]{Remark}
\newcommand{\bdfn}{\begin{defn}}
\newcommand{\bthm}{\begin{thm}}
\newcommand{\blmma}{\begin{lemma}}
\newcommand{\bppsn}{\begin{prop}}
\newcommand{\bcrlre}{\begin{corr}}
\newcommand{\bxmpl}{\begin{xmpl}}
\newcommand{\brmrk}{\begin{rmk}}
\newcommand{\edfn}{\end{defn}}
\newcommand{\ethm}{\end{thm}}
\newcommand{\elmma}{\end{lemma}}
\newcommand{\eppsn}{\end{prop}}
\newcommand{\ecrlre}{\end{corr}}
\newcommand{\exmpl}{\end{xmpl}}
\newcommand{\ermrk}{\end{rmk}}
\newcommand{\IC}{\mathbb{C}}

\newcommand{\homc}{\mathrm{Hom}_\IC}
\newcommand{\id}{\mathrm{id}}

\newcommand{\twoform}{{{\Omega}^2}( \mathcal{A} )}
\newcommand{\tensora}{\otimes_{\mathcal{A}}}

\newcommand{\tensorc}{\otimes_{\mathbb{C}}}
\newcommand{\A}{\mathcal{A}}

\newcommand{\E}{\mathcal{E}}

\newcommand{\F}{\mathcal{F}}

\newcommand{\zeroE}{{}_0\E}

\newcommand{\Psym}{P_{\rm sym}}
\newcommand{\Hom}{{\rm Hom}}

\newcommand{\RNum}[1]{\uppercase\expandafter{\romannumeral #1\relax}}

\title{Levi-Civita connection for $SU_q(2)$}
\author{Sugato Mukhopadhyay}
\address{\normalfont Indian Statistical Institute\\
	203 B.T. Road, 
	Kolkata, India
}
\email{m.xugato@gmail.com}

\begin{document}

\begin{abstract}
We prove that the $4D_\pm$ calculi on the quantum group $SU_q(2)$ satisfy a metric-independent sufficient condition for the existence of a unique bicovariant Levi-Civita connection corresponding to every bi-invariant pseudo-Riemannian metric.
\end{abstract}

\maketitle

\section{Introduction}
The quantum group $SU_q(2)$ introduced in \cite{woronowicz2} and the notion of bicovariant differential calculi on Hopf algebras was introduced in \cite{woronowicz} by Woronowicz.
The question of bicovariant Levi-Civita connections on bicovariant differential calculi of compact quantum groups have been investigated by Heckenberger and Schm{\"u}dgen in \cite{heckenberger} for the quantum groups $SL_q(N)$, $O_q(N)$ and $Sp_q(N)$. On the other hand, Beggs, Majid and their collaborators have studied Levi-Civita connections on quantum groups and homogeneous spaces in various articles, and a comprehensive account can be found in \cite{beggsmajidbook}.

More recently, in \cite{article6}, bicovariant connections on arbitrary bicovariant differential calculi of compact quantum groups and the notion of their metric compatibility with respect to arbitrary bi-invariant pseudo-Riemannian metrics was studied. In that article, the construction of a canonical bicovariant torsionless connection on a calculus was presented (Theorem 5.3 of \cite{article6}), provided that Woronowicz's braiding map $\sigma$ for the calculus satisfies a diagonalisability condition. Also, a metric-independent sufficient condition for the existence of a unique bicovariant Levi-Civita connection (in the sense of Definition 6.3 of \cite{article6}) was provided in Theorem 7.9 of \cite{article6}.

In this article, we will investigate the theory of \cite{article6} in the context of the $4D_\pm$ calculi of the compact quantum group $SU_q(2)$ which were explicitly described in \cite{woronowicz2} and then \cite{stachura}. In Section \ref{section0}, we recall the notion of covariant Levi-Civita connections on bicovariant differential calculi as formulated in \cite{article6}. In Section \ref{section2}, the $4D_\pm$ calculi on $SU_q(2)$ are recalled and we show that Woronowicz's braiding map for the $4D_\pm$ calculi satisfy the diagonalisability condition mentioned above. In Section \ref{section4}, we construct a bicovariant torsionless connection. In Section \ref{section5}, we will show that the sufficiency condition of Theorem 7.9 of \cite{article6} is satisfied by both calculi, except for at most finitely many values of $q$, and hence we can conclude the existence of a unique bicovariant Levi-Civita connection, corresponding to a bi-invariant pseudo-Riemannian metric.

\section{Levi-Civita connections on bicovariant differential calculi} \label{section0}

In this section, we recall the notion of Levi-Civita connections on bicovariant differential calculi as formulated in \cite{article6}.

We say that $(\E, \Delta_\E, {}_\E \Delta)$ is a bicovariant bimodule over a Hopf algebra $\A$ if $\E$ is a bimodule over $\A$, $(\E, \Delta_\E)$ is a left $\A$-comodule, $(\E, {}_\E \Delta)$ is a right $\A$-comodule, subject to the following compatibility conditions:
\begin{eqnarray*}
	\Delta_\E(a \rho) = \Delta(a) \Delta_\E(\rho), & \Delta_\E( \rho a) =  \Delta_\E(\rho) \Delta(a)\\
	{}_\E \Delta(a \rho) = \Delta(a) {}_\E \Delta(\rho), & {}_\E \Delta(\rho a) = {}_\E \Delta(\rho) \Delta(a),
\end{eqnarray*}
where $\rho$ is an arbitrary element of $\E$ and $a$ is an arbitrary element of $\A$. If $(\E, \Delta_\E, {}_\E \Delta)$ is a bicovariant bimodule over $\A$, we say that an element $e$ in $\E$ is left (respectively, right) invariant if $\Delta_\E(e) = 1 \tensorc e$ (respectively, ${}_\E \Delta(e) = e \tensorc 1$). In this article, we will denote the vector space of elements of $\E$ invariant under the left coaction of $\A$ by ${}_0 \E$, and that of elements invariant under the right coaction of $\A$ by $\E_0$.
If $\E$ and $\F$ are two bicovariant bimodules over $\A$, a $\IC$-linear map $T: \E \to \F$ is said to be left covariant if $\Delta_\F \circ T = (\id \tensorc T) \circ \Delta_\E$. $T$ is said to be right covariant if ${}_\F \Delta \circ T = (T \tensorc \id) \circ {}_\E \Delta$. $T$ is called bicovariant if it is both left-covariant and right-covariant.

A (first order) differential calculus $(\E, d)$ over a Hopf algebra $\A$ is called a bicovariant differential calculus if the following conditions are satisfied:
	\begin{itemize}
		\item[(i)] For any $a_k,b_k$ in $\A$, $k=1, \dots, K$,
		$$(\sum_k a_kdb_k = 0) \text{ implies that }(\sum_k \Delta(a_k)(\id \tensorc d)\Delta(b_k)=0),$$
		\item[(ii)] For any $a_k,b_k$ in $\A$, $k=1, \dots, K$,
		$$(\sum_k a_kdb_k = 0) \text{ implies that } (\sum_k \Delta(a_k)(d \tensorc \id)\Delta(b_k)=0).$$
	\end{itemize}

Woronowicz (\cite{woronowicz}) proved that a bicovariant differential calculus is endowed with canonical left and right-comodule coactions of $\A$, making it into a bicovariant bimodule $(\E, \Delta_\E, {}_\E \Delta)$. Moreover, the map $d: \A \to \E$ is a bicovariant map.

Next, we state the construction of the associated space of two-forms, $\Omega^2(\A)$ for a bicovariant differential calculus of an arbitrary unital Hopf algebra $\A$ as in \cite{woronowicz}. To do so, we need to recall the braiding map $\sigma$ for bicovariant bimodules.

\begin{prop}{(Proposition 3.1 of \cite{woronowicz})} \label{4thmay20193}
	Given a bicovariant bimodule $\E$ on a Hopf algebra $\A$, there exists a unique bimodule homomorphism
	$$\sigma: \E \tensora \E \to \E \tensora \E ~ {\rm such} ~ {\rm that} $$ 
	\begin{equation} \label{30thapril20191} \sigma(\omega \tensora \eta)= \eta \tensora \omega \end{equation}
	for any left-invariant element $\omega$ and right-invariant element $\eta$ in $\E$, under the coactions of $\A$. $\sigma$ is an invertible bicovariant $\A$-bimodule map from $\E \tensora \E$ to itself. Moreover, $\sigma$ satisfies the following braid equation on $\E \tensora \E \tensora \E:$
	$$ (\id \tensora \sigma)(\sigma \tensora \id)(\id \tensora \sigma)= (\sigma \tensora \id)(\id \tensora \sigma)(\sigma \tensora \id). $$
\end{prop}

Given a bicovariant first order differential calculus $(\E, d)$, there exists a unique braiding map $ \sigma $, by Proposition \ref{4thmay20193}. The space of two-forms is defined to be the bicovariant bimodule
	$$\twoform := (\E \tensora \E) \big/ \rm {\rm Ker} (\sigma - 1).$$
	The symbol $\wedge$ denotes the quotient map, which is a bicovariant bimodule map,
	\[ \wedge: \E\tensora \E \to  \twoform. \]
	The map $d: \A \to \E$ extends to a unique exterior derivative map (to be denoted again by $d$),
	\[d: \E \to \twoform, \]
	such that, for all $a$ in $\A$ and $\rho$ in $\E$,
	\begin{itemize}
		\item[(i)] $d(a \rho) = da \wedge \rho + a d (\rho),$
		\item[(ii)] $d(\rho a) = d (\rho) a - \rho \wedge da,$
		\item[(iii)] $d$ is bicovariant.
	\end{itemize}

Let us, from now on, denote the subspace of left-invariant elements of an arbitrary bicovariant bimodule $\E$ by the symbol ${}_0 \E$. By Proposition 2.5 of \cite{article5}, the vector space ${}_0 \E \tensorc {}_0 \E$ can be identified with the space ${}_0(\E \tensora \E)$ of left-invariant elements of $\E\tensora \E$. The isomorphism ${}_0 \E \tensorc {}_0 \E \to {}_0 (\E \tensora \E)$ is given by
\begin{equation} \label{24thjan20192} \omega_i \tensorc \omega_j \mapsto \omega_i \tensora \omega_j \end{equation},
where $\{  \omega_i \}_i$ is a vector space basis ${}_0 \E$.\\
Moreover, by the bicovariance of the map $\sigma: \E \tensora \E \to \E \tensora \E$, we get the restriction (see Equation 18 of \cite{article6}):
\begin{equation} \label{21staugust20195} {}_0 \sigma := \sigma|_{{}_0 \E \tensorc {}_0 \E} : {}_0 \E \tensorc {}_0 \E \to {}_0 \E \tensorc {}_0 \E. \end{equation}
From now on, we are going to work under the assumption that ${}_0 \sigma$ is a diagonalisable map between finite dimensional vector spaces. In \cite{article6}, this assumption was crucially used to set up a framework for the existence of a unique bicovariant Levi-Civita connection on a bicovariant differential calculus satisfying the assumption.

Let us introduce some notations and definitions so that we can recall the framework mentioned above.
\begin{defn} \label{22ndaugust20192}
	Suppose the map ${}_0 \sigma$ is diagonalisable. The eigenspace decomposition of ${}_0 \E \tensorc {}_0 \E$ will be denoted by ${}_0 \E \tensorc {}_0 \E = \bigoplus_{\lambda \in \Lambda} V_\lambda$, where $\Lambda$ is the set of distinct eigenvalues of ${}_0 \sigma$ and $V_\lambda$ is the eigenspace of ${}_0\sigma$ corresponding to the eigenvalue $\lambda$. Thus, for example, $ V_1 $ will denote the eigenspace of $ {}_0 \sigma $ for the eigenvalue $ \lambda = 1. $
	
	Moreover,	we define ${}_0 \E \tensorc^{\rm sym} {}_0 \E $ to be the eigenspace of ${}_0 \sigma$ with eigenvalue $1$, i.e.,
	\[ {}_0 \E \tensorc^{\rm sym} {}_0 \E := V_1.\]
	We also define ${}_0 \F := \bigoplus_{\lambda \in \Lambda \backslash \{1\}} V_\lambda.$
	Finally, we will denote by ${}_0(\Psym)$ the idempotent element in Hom$({}_0 \E \tensorc {}_0 \E, {}_0 \E \tensorc {}_0 \E)$ with range ${}_0 \E \tensorc^{\rm sym} {}_0 \E$ and kernel ${}_0\F$.
\end{defn}

Lets us introduce the notion of bi-invariant pseudo-Riemannian metric on a bicovariant $\A$-bimodule $\E$.

\begin{defn}{(\cite{heckenberger}, Definition 4.1 of \cite{article6})} \label{24thmay20191}
	Suppose $ \E $ is a bicovariant $\A$ bimodule and $ \sigma: \E \tensora \E \rightarrow \E \tensora \E $ be the map as in Proposition \ref{4thmay20193}. A bi-invariant pseudo-Riemannian metric for the pair $ (\E, \sigma) $ is a right $\A$-linear map $g:\E \tensora \E \to \A$ such that the following conditions hold:
	\begin{enumerate}
		\item[(i)] $ g \circ \sigma = g. $
		\item[(ii)] If $g(\rho \tensora \nu)=0$ for all $\nu$ in $\E,$ then $\rho = 0.$ 
		\item[(iii)] The map $g$ is bi-invariant, i.e for all $\rho, \nu$ in $\E$,
		\[ (\id \tensorc \epsilon g)(\Delta_{(\E \tensora \E)}(\rho \tensora \nu)) = g(\rho \tensora \nu), \]
		\[ (\epsilon g \tensorc \id)({}_{(\E \tensora \E)}\Delta(\rho \tensora \nu)) = g(\rho \tensora \nu). \]
	\end{enumerate}
\end{defn}

Now we can define the torsion of a connection and the compatibility of a left-covariant connection with a bi-invariant pseudo-Riemannian metric.

\begin{defn} (\cite{heckenberger})
	Let $ (\E, d)$ be a bicovariant differential calculus on $\A.$ A (right) connection on $\E$ is a $\IC$-linear map $\nabla: \E \to \E \tensora \E$ such that, for all $a$ in $\A$ and $\rho$ in $\E$, the following equation holds:
	$$ \nabla(\rho a)= \nabla(\rho)a + \rho \tensora da. $$
	The map $\nabla$ is said to a left-covariant, right-covariant or bicovariant connection is it is a left-covariant, right-covariant or bicovariant map, respectively.
	The torsion of a connection $\nabla$ on $\E$ is the right $\A$-linear map 
	$$T_\nabla:= \wedge \circ \nabla + d :\E \to \twoform.$$
	The connection $\nabla$ is said to be torsionless if $T_{\nabla}=0$.
\end{defn}

Our notion of torsion is the same as that of \cite{heckenberger}, with the only difference being that they work with left connections.

\begin{defn}{(Definitions 6.1 and 6.3 of \cite{article6})}
	Let $\nabla$ be a left-covariant connection on a bicovariant calculus $(\E, d)$ such that the map ${}_0 \sigma$ is diagonalisable, and $ g $ a bi-invariant pseudo-Riemannian metric. Then we define 
	$$\widetilde{\Pi^0_g}(\nabla):{}_0\E \tensorc {}_0\E \rightarrow {}_0\E ~ {\rm by} ~ {\rm the} ~ {\rm following} ~ {\rm formula}:$$
	\begin{equation} \label{1stmay20191}	\widetilde{\Pi^0_g}(\nabla)(\omega_i \tensorc \omega_j)=2(\id \tensorc g)(\sigma \tensorc \id)(\nabla \tensorc \id){}_0(\Psym)(\omega_i \tensorc \omega_j). \end{equation}
	Next, for all $ \omega_1,\omega_2$ in $\zeroE $ and $a$ in $\A$, we define $\widetilde{\Pi_g}(\nabla):\E \tensora \E \rightarrow \E$ by 
	$$\widetilde{\Pi_g}(\nabla) (\omega_1 \tensora \omega_2 a)= \widetilde{\Pi^0_g}(\nabla)(\omega_1 \tensorc \omega_2)a + g(\omega_1 \tensora \omega_2)da.$$
	Finally, $\nabla$ is said to be compatible with $g$, if, as maps from $\E \tensora \E$ to $\E$,
	\[ \widetilde{\Pi_g}(\nabla) = dg. \]
\end{defn}

This allows us to give the definition of a Levi-Civita connection.

\begin{defn}
	Let $ (\E, d) $ be a bicovariant differential calculus such that the map $ {}_0 \sigma $ is diagonalisable and $ g $ a pseudo-Riemannian bi-invariant metric on $ \E. $ A left-covariant connection $ \nabla $ on $ \E $ is called a Levi-Civita connection for the triple $ (\E, d, g) $ if it is torsionless and compatible with $g.$ 
\end{defn}

In \cite{article6}, it was shown that this suitably generalises the notion of Levi-Civita connections for bicovariant differential calculi on Hopf algebras.

Then, we have the following metric-independent sufficient condition for the existence of a unique bicovariant Levi-Civita connection.

\begin{thm}{(Theorem 7.9 of \cite{article6})} \label{15thjune20193}
	Suppose $(\E, d)$ is a bicovariant differential calculus over a cosemisimple Hopf algebra $\A$ such that the map ${}_0 \sigma$ is diagonalisable and $ g $ be a bi-invariant pseudo-Riemannian metric. If the map \[ ({}_0(\Psym))_{23}: (\zeroE \tensorc^{\rm sym} \zeroE) \tensorc \zeroE \rightarrow \zeroE \tensorc (\zeroE \tensorc^{\rm sym} \zeroE)  \] is an isomorphism, then the triple $ (\E, d, g) $ admits a unique bicovariant Levi-Civita connection.
\end{thm}

\section{The $4D_\pm$ calculi on $SU_q(2)$ and the braiding map} \label{section2}

In this section we recall briefly the definition quantum group $SU_q(2)$ and the $4D_\pm$ calculi on $SU_q(2)$. Then we show that the map ${}_0 \sigma : {}_0\ E \tensorc {}_0 \E \to {}_0\ E \tensorc {}_0 \E$ is actually diagonalisable. Our main reference for the details is \cite{stachura}.

For $q \in [-1,1] \backslash {0}$, $SU_q(2)$ is the C*-algebra generated by the two elements $\alpha$ and $\gamma$, and their adjoints, satisfying the following relations:
\begin{align*}
	\alpha^* \alpha + \gamma^* \gamma = 1 , \quad
	\alpha \alpha^* + q^2 \gamma \gamma^* = 1,\\
	\gamma^* \gamma = \gamma \gamma^*, \quad
	\alpha \gamma = q \gamma \alpha, \quad
	\alpha \gamma^* = q \gamma^* \alpha.
\end{align*} 
The comultiplication map $\Delta$ is given by
\[ \Delta(\alpha) = \alpha \tensorc \alpha - q \gamma^\ast \tensorc \gamma, \quad \Delta(\gamma) = \gamma \tensorc \alpha + \alpha^\ast \tensorc \gamma. \]
This makes $SU_q(2)$ into a compact quantum group. We will denote the Hopf $*$-algebra generated by the elements $\alpha$, $\gamma$ by the symbol $\A$.\\
 
In \cite{stachura}, it is explicitly proven that there does not exists any three-dimensional bicovariant differential calculi and exactly two inequivalent four-dimensional calculi for $SU_q(2)$. We use the description of the two bicovariant calculi, $4 {\rm D}_+$ and $4 {\rm D}_-$, as given in \cite{stachura}. We will rephrase some of the notation to fit our formalism. For $q \in (-1,1) \backslash \{0\}$, the first order differential calculi $\E$ of both the $4 {\rm D}_+$ and $4 {\rm D}_-$ calculi are  bicovariant $\A$-bimodules such that the space ${}_0 \E$ of one-forms invariant under the left coaction of $\A$ is a $4$-dimensional vector space. We will denote a preferred basis of ${}_0 \E$ by $\{ \omega_i \}_{i =1,2,3,4}$. Here we have replaced the notation in \cite{stachura} with $\omega_i = \Omega_i$.

The following is the explicit description of the exterior derivative $d$ on ${}_0\E$ for the preferred basis $\{ \omega_i \}_{i=1}^4$ mentioned above.
\begin{prop}{(Equation (5.2) of \cite{stachura})} \label{24thjan20191}
	Let $d: \E \to \twoform$ be the exterior derivative of the $4D_\pm$ calculus.
	\begin{align*}
		d(\omega_1) &= \pm \sqrt{r}\omega_1 \wedge \omega_3, & d(\omega_2) =& \mp \frac{\sqrt{r}}{q^2} \omega_2 \wedge \omega_3,\\
		d(\omega_3) &= \pm \frac{\sqrt{r}}{q} \omega_1 \wedge \omega_2, & d(\omega_4) =& 0,
	\end{align*}
where the upper sign stand for $4D_{+}$ and the lower for $4D_{-}$, and $r = 1 + q^2$.
\end{prop}

Now we show that the map ${}_0\sigma$ for $SU_q(2)$ satisfies the diagonalisability condition by giving explicit bases for eigenspaces of ${}_0 \sigma$. We will use the explicit action of $\sigma$ on elements $\omega_i \tensora \omega_j$, $i,j =1,2,3,4$ as given in Equation (4.1) of \cite{stachura}

\begin{prop} \label{26thjan20201}
	For $SU_q(2)$, the map ${}_0 \sigma$ is diagonalisable and has the minimal polynomial equation
	\[ ({}_0 \sigma - 1)({}_0 \sigma + q^2)({}_0 \sigma + q^{-2})=0. \]
\end{prop}
\begin{proof}
	The proof of this result is by explicit listing of eigenvectors of ${}_0 \sigma$ for eigenvalues $1, q^2, q^{-2}$ and by a dimension argument. Throughout we make use of the canonical identification $\omega_i \tensorc \omega_j \mapsto \omega_i \tensora \omega_j$ as stated in \eqref{24thjan20192}.\\
	Either by directly applying ${}_0 \sigma$ on the following linearly independent two-tensors or from Equation (4.2) of \cite{stachura}, we get that the following are in the eigenspace corresponding to eigenvalue $1$:
	\begin{eqnarray*}
	&\omega_1 \tensorc \omega_1, \omega_2 \tensorc \omega_2, \omega_3 \tensorc \omega_3 + t \omega_1 \tensorc \omega_2, \omega_4 \tensorc \omega_4,\\
	&\omega_1 \tensorc \omega_2 + \omega_2 \tensorc \omega_1, \omega_2 \tensorc \omega_3 + q^2 \omega_3 \tensorc \omega_2,\\
	&q^2 \omega_1 \tensorc \omega_3 + \omega_3 \tensorc \omega_1, \frac{t^2 k}{q^2 \sqrt{r}} \omega_2 \tensorc \omega_3 - \omega_2 \tensorc \omega_4 - \omega_4 \tensorc \omega_2,\\
	&\frac{t^2 k}{\sqrt{r}} \omega_1 \tensorc \omega_3 + \omega_1 \tensorc \omega_4 + \omega_4 \tensorc \omega_1, \frac{t^2 k}{q \sqrt{r}} \omega_1 \tensorc \omega_2 + \omega_3 \tensorc \omega_4 + \omega_4 \tensorc \omega_3.
	\end{eqnarray*}
	By explicit computation, the following linearly independent two-tensors are in the eigenspace corresponding to the eigenvalue $q^2$:
	\begin{eqnarray*}
		& \frac{tqk}{\sqrt{r}} \omega_2 \tensorc \omega_3 - q^2 \omega_2 \tensorc \omega_4 - \frac{tk}{q \sqrt{r}} \omega_3 \tensorc \omega_2 + \omega_4 \tensorc \omega_2,\\
		& -\frac{tk}{q \sqrt{r}} \omega_1 \tensorc \omega_3 - q^2 \omega_1 \tensorc \omega_4 + \frac{tqk}{\sqrt{r}} \omega_3 \tensorc \omega_1 + \omega_4 \tensorc \omega_1,\\
		&- \frac{tk}{\sqrt{r}} \omega_1 \tensorc \omega_2 + \frac{tk}{\sqrt{r}} \omega_2 \tensorc \omega_1 + \frac{t^2 k}{\sqrt{r}}\omega_3 \tensorc \omega_3 - q^2 \omega_3 \tensorc \omega_4 + \omega_4 \tensorc \omega_3.		
	\end{eqnarray*}
	By explicit computation, the following linearly independent two-tensors are in the eigenspace corresponding to the eigenvalue $q^{-2}$:
	\begin{eqnarray*}
		& \frac{tqk}{\sqrt{r}} \omega_2 \tensorc \omega_3 + \omega_2 \tensorc \omega_4 - \frac{tk}{q \sqrt{r}} \omega_3 \tensorc \omega_2 - q^{-2} \omega_4 \tensorc \omega_2,\\
		& -\frac{tk}{q \sqrt{r}} \omega_1 \tensorc \omega_3 + \omega_1 \tensorc \omega_4 + \frac{tqk}{\sqrt{r}} \omega_3 \tensorc \omega_1 - q^{-2} \omega_4 \tensorc \omega_1,\\
		&- \frac{tk}{\sqrt{r}} \omega_1 \tensorc \omega_2 + \frac{tk}{\sqrt{r}} \omega_2 \tensorc \omega_1 + \frac{t^2 k}{\sqrt{r}}\omega_3 \tensorc \omega_3 + \omega_3 \tensorc \omega_4 - q^{2} \omega_4 \tensorc \omega_3.		
	\end{eqnarray*}
	We have thus accounted for $16$ linearly independent elements of ${}_0 \E \tensorc {}_0 \E$. Since ${}_0 \E$ has dimension $4$, ${}_0 \E \tensorc {}_0 \E$ has dimension $16$. Hence we have a basis, and in particular bases for the eigenspace decomposition, of ${}_0 \E \tensorc {}_0 \E$. Moreover, ${}_0 \sigma$ satisfies the minimal polynomial
	\[ ({}_0 \sigma - 1)({}_0 \sigma + q^2)({}_0 \sigma + q^{-2})=0. \]
\end{proof}

\section{A bicovariant torsionless connection} \label{section4}

In this section, using the fact that ${}_0 \sigma$ is diagonalisable and ${}_0\E \tensorc {}_0 \E$ admits an eigenspace decomposition, we construct a bicovariant torsionless connection on the $4D_\pm$ calculus.

\begin{rmk} \label{26thjan20202}
	Note that since any element $\rho$ in $\E$ can be uniquely expressed as $\rho = \sum_i \omega_i a_i$ for some $a_i$ in $\A$ (Theorem 2.1 of \cite{woronowicz}), a connection on $\E$ is determined by its action on the basis $\{\omega_i\}_i$.
\end{rmk}

By Proposition \ref{26thjan20201}, we have the eigenspace decomposition
\[ {}_0 \E \tensorc {}_0 \E = {\rm Ker}({}_0 \sigma - \id) \oplus {\rm Ker}({}_0 \sigma + q^2) \oplus {\rm Ker}({}_0 \sigma + q^{-2}). \]
Since ${\rm Ker}(\wedge) = {\rm Ker}({}_0 \sigma - \id)$, we have that
\[ {\rm Ker}({}_0 \sigma + q^2) \oplus {\rm Ker}({}_0 \sigma + q^{-2})  \cong \twoform , \]
with the isomorphism being given by $\wedge|_{{\rm Ker}({}_0 \sigma + q^2) \oplus {\rm Ker}({}_0 \sigma + q^{-2})}$. Let us denote ${\rm Ker}({}_0 \sigma + q^2) \oplus {\rm Ker}({}_0 \sigma + q^{-2})$ by ${}_0 \F$ from now on. This is consistent with the notation adopted in Definition \ref{22ndaugust20192}.

\begin{thm} \label{31stmar20201}
	Let $\{ \omega_i \}_i$ be the preferred basis for the $4D_\pm$ calculus on $SU_q(2)$. For $i = 1,2,3,4$, we define 
	\[ \nabla_0(\omega_i) = - (\wedge|_{{}_0 \F})^{-1} \circ d(\omega_i) \in {}_0 \E \tensorc {}_0 \E. \]
	Then, $\nabla$ extends to a bicovariant torsionless connection on $\E$. More explicitly,
	\begin{equation*}
		\begin{aligned}
		&\nabla_0(\omega_1) &=& \mp \frac{{r}}{t^2 k (q^2 + 1)^2}\big(\frac{2tk}{q \sqrt{r}} \omega_1 \tensorc \omega_3 + tq \omega_1 \tensorc \omega_4 - \frac{2tqk}{\sqrt{r}} \omega_3 \tensorc \omega_1 + tq \omega_4 \tensorc \omega_1\big) \\
		&\nabla_0(\omega_2) &=& \pm \frac{{r}}{t^2 k (q^2 + 1)^2}\big( \frac{2tqk}{ \sqrt{r}} \omega_2 \tensorc \omega_3 - tq \omega_2 \tensorc \omega_4 - \frac{2tk}{q\sqrt{r}} \omega_3 \tensorc \omega_2 - tq \omega_4 \tensorc \omega_2 \big) \\
		&\nabla_0(\omega_3) &=& \pm \frac{q {r}}{tk (q^2 + 1)^2}\big( \frac{2tk}{\sqrt{r}} \omega_1 \tensorc \omega_2 - \frac{2tk}{\sqrt{r}} \omega_2 \tensorc \omega_1 - \frac{t^2 k}{\sqrt{r}}\omega_3 \tensorc \omega_3\\ &&&+ tq \omega_3 \tensorc \omega_4 + tq \omega_4 \tensorc \omega_3 \big) \\
		&\nabla_0(\omega_4) &=& \ 0
		\end{aligned}
	\end{equation*}
\end{thm}
\begin{proof}
	By the definition of $\nabla_0$,
	\[ \wedge \circ \nabla_0(\omega_i) = - \wedge \circ (\wedge|_{{}_0 \F})^{-1} \circ d(\omega_i) = - d (\omega_i). \]
	Therefore, for any element $\rho = \sum_i \omega_i a_i$ in $\E$,
	\begin{eqnarray*}
	&&\wedge \circ \nabla_0(\sum_i \omega_i a_i)  = \wedge \circ \sum_i(\nabla_0( \omega_i) a_i + \omega_i \tensora a_i )  \\
	= && - \sum_i\big(\wedge \circ (\wedge|_{{}_0 \F})^{-1} \circ d(\omega_i) a_i + \omega_i \wedge a_i \big) \\
	 = && - \sum_i( d(\omega_i) a_i + \omega_i \wedge a_i) = - \sum_i d(\omega_i a_i).
	\end{eqnarray*}
	Hence $\nabla_0$ is a torsionless connection. The construction of $\nabla_0$ is the same as that in Theorem 5.3 of \cite{article6}. Hence, by that theorem, our connection $\nabla_0$ is bicovariant.
	
	Now we derive $\nabla_0$ explicitly on each $\omega_i$ using the formulas for $d(\omega_i)$ in Proposition \ref{24thjan20191}.\\
	We have that $d(\omega_1) = \pm \sqrt{r} \omega_1 \wedge \omega_3$. The decomposition of $\omega_1 \tensorc \omega_3$ as a linear combination of the basis eigenvectors listed in Proposition \ref{26thjan20201} is given by
	\begin{equation*}
	\begin{aligned}
		\omega_1 \tensorc \omega_3 =& \ \frac{2q^2}{(q^2 + 1)^2}\big(q^2 \omega_1 \tensorc \omega_3 + \omega_3 \tensorc \omega_1 \big)\\
		 &- \frac{q^2 \sqrt{r}}{k(q^2 + 1)^2}\big(\frac{t^2 k}{\sqrt{r}} \omega_1 \tensorc \omega_3 + \omega_1 \tensorc \omega_4 + \omega_4 \tensorc \omega_1\big)\\
		&-\frac{\sqrt{r}}{t^2 k (q^2 + 1)^2}\big(-\frac{tk}{q \sqrt{r}} \omega_1 \tensorc \omega_3 - q^2 \omega_1 \tensorc \omega_4 + \frac{tqk}{\sqrt{r}} \omega_3 \tensorc \omega_1 + \omega_4 \tensorc \omega_1\big)\\
		&-\frac{\sqrt{r}}{t^2 k (q^2 + 1)^2}\big(-\frac{tk}{q \sqrt{r}} \omega_1 \tensorc \omega_3 + \omega_1 \tensorc \omega_4 + \frac{tqk}{\sqrt{r}} \omega_3 \tensorc \omega_1 - q^2 \omega_4 \tensorc \omega_1\big).
	\end{aligned}
	\end{equation*}
	Since the first two terms in the above decomposition are elements of ${\rm Ker}({}_0 \sigma - \id) = {\rm Ker}(\wedge)$, applying $\wedge$ on both sides, we have
	\begin{equation*}
	\begin{aligned}
		\omega_1 \wedge \omega_3 =& \wedge \Big(-\frac{\sqrt{r}}{t^2 k (q^2 + 1)^2}\big(-\frac{tk}{q \sqrt{r}} \omega_1 \tensorc \omega_3 - q^2 \omega_1 \tensorc \omega_4 + \frac{tqk}{\sqrt{r}} \omega_3 \tensorc \omega_1 + \omega_4 \tensorc \omega_1\big)\\
		&-\frac{\sqrt{r}}{t^2 k (q^2 + 1)^2}\big(-\frac{tk}{q \sqrt{r}} \omega_1 \tensorc \omega_3 + \omega_1 \tensorc \omega_4 + \frac{tqk}{\sqrt{r}} \omega_3 \tensorc \omega_1 - q^2 \omega_4 \tensorc \omega_1\big)\Big),
	\end{aligned}
	\end{equation*}
	and since the last two terms in the decomposition are from ${}_0 \F$,
	\begin{equation*}
	\begin{aligned}
		(\wedge|_{{}_0 \F})^{-1}(\omega_1 \wedge \omega_3) = &-\frac{\sqrt{r}}{t^2 k (q^2 + 1)^2}\big(-\frac{tk}{q \sqrt{r}} \omega_1 \tensorc \omega_3 - q^2 \omega_1 \tensorc \omega_4 + \frac{tqk}{\sqrt{r}} \omega_3 \tensorc \omega_1 + \omega_4 \tensorc \omega_1\big)\\
		&-\frac{\sqrt{r}}{t^2 k (q^2 + 1)^2}\big(-\frac{tk}{q \sqrt{r}} \omega_1 \tensorc \omega_3 + \omega_1 \tensorc \omega_4 + \frac{tqk}{\sqrt{r}} \omega_3 \tensorc \omega_1 - q^2 \omega_4 \tensorc \omega_1\big).
	\end{aligned}
	\end{equation*}
	Thus, by the construction of $\nabla_0$, we have
	\begin{equation*}
	\begin{aligned}
		\nabla_0(\omega_1) =& \mp  \big( -\frac{{r}}{t^2 k (q^2 + 1)^2}\big(-\frac{tk}{q \sqrt{r}} \omega_1 \tensorc \omega_3 - q^2 \omega_1 \tensorc \omega_4 + \frac{tqk}{\sqrt{r}} \omega_3 \tensorc \omega_1 + \omega_4 \tensorc \omega_1\big)\\
		&-\frac{{r}}{t^2 k (q^2 + 1)^2}\big(-\frac{tk}{q \sqrt{r}} \omega_1 \tensorc \omega_3 + \omega_1 \tensorc \omega_4 + \frac{tqk}{\sqrt{r}} \omega_3 \tensorc \omega_1 - q^2 \omega_4 \tensorc \omega_1\big) \big)\\
		=& \mp \frac{{r}}{t^2 k (q^2 + 1)^2}\big(\frac{2tk}{q \sqrt{r}} \omega_1 \tensorc \omega_3 + tq \omega_1 \tensorc \omega_4 - \frac{2tqk}{\sqrt{r}} \omega_3 \tensorc \omega_1 + tq \omega_4 \tensorc \omega_1\big)
	\end{aligned}
	\end{equation*}
Proposition \ref{24thjan20191} also gives that $ d(\omega_2) = \mp \frac{\sqrt{r}}{q^2} \omega_2 \wedge \omega_3,
d(\omega_3) = \pm \frac{\sqrt{r}}{q} \omega_1 \wedge \omega_2$ and $ d(\omega_4) = 0 $. So, similarly, we have
\begin{equation*}
	\begin{aligned}
		\omega_2 \tensorc \omega_3 = & \frac{2}{(q^2 + 1)^2} \big( \omega_2 \tensorc \omega_3 + q^2 \omega_3 \tensorc \omega_2 \big) - \frac{q^4 \sqrt{r}}{k(q^2 + 1)^2}\big(\frac{t^2 k}{\sqrt{r}} \omega_2 \tensorc \omega_3 - \omega_2 \tensorc \omega_4 - \omega_4 \tensorc \omega_2 \big)\\
		&+\frac{q^2 \sqrt{r}}{t^2 k (q^2 + 1)^2}\big(\frac{tqk}{ \sqrt{r}} \omega_2 \tensorc \omega_3 - q^2 \omega_2 \tensorc \omega_4 - \frac{tk}{q\sqrt{r}} \omega_3 \tensorc \omega_2 + \omega_4 \tensorc \omega_2 \big)\\
		&+\frac{q^2 \sqrt{r}}{t^2 k (q^2 + 1)^2}\big(\frac{tqk}{ \sqrt{r}} \omega_2 \tensorc \omega_3 + \omega_2 \tensorc \omega_4 - \frac{tk}{q \sqrt{r}} \omega_3 \tensorc \omega_2 - q^2 \omega_4 \tensorc \omega_1\big),
	\end{aligned}
\end{equation*}
and hence,
\begin{equation*}
	\begin{aligned}
		\nabla_0(\omega_2) = \pm \frac{{r}}{t^2 k (q^2 + 1)^2}\big( \frac{2tqk}{ \sqrt{r}} \omega_2 \tensorc \omega_3 - tq \omega_2 \tensorc \omega_4 - \frac{2tk}{q\sqrt{r}} \omega_3 \tensorc \omega_2 - tq \omega_4 \tensorc \omega_2 \big).
	\end{aligned}
\end{equation*}
Moreover,
\begin{equation*}
	\begin{aligned}
		\omega_1 \tensorc \omega_2 = & \frac{2tq^2}{(q^2 + 1)^2}\big( \omega_1 \tensorc \omega_2 + \omega_2 \tensorc \omega_1 \big) + \frac{2q^2}{(q^2 + 1)^2}\big(\omega_3 \tensorc \omega_3 + t \omega_1 \tensorc \omega_2\big)\\
		&- \frac{q^3 \sqrt{r}}{k(q^2 + 1)^2}\big( \frac{t^2 k}{q \sqrt{r}} \omega_1 \tensorc \omega_2 + \omega_3 \tensorc \omega_4 + \omega_4 \tensorc \omega_3 \big)\\
		& - \frac{q^2 \sqrt{r}}{tk (q^2 + 1)^2}\big(- \frac{tk}{\sqrt{r}} \omega_1 \tensorc \omega_2 + \frac{tk}{\sqrt{r}} \omega_2 \tensorc \omega_1 + \frac{t^2 k}{\sqrt{r}}\omega_3 \tensorc \omega_3 - q^2 \omega_3 \tensorc \omega_4 + \omega_4 \tensorc \omega_3 \big)\\
		& - \frac{q^2 \sqrt{r}}{tk (q^2 + 1)^2}\big( - \frac{tk}{\sqrt{r}} \omega_1 \tensorc \omega_2 + \frac{tk}{\sqrt{r}} \omega_2 \tensorc \omega_1 + \frac{t^2 k}{\sqrt{r}}\omega_3 \tensorc \omega_3 + \omega_3 \tensorc \omega_4 - q^{2} \omega_4 \tensorc \omega_3 \big),
	\end{aligned}
\end{equation*}
and hence,
\begin{equation*}
	\begin{aligned}
		\nabla_0(\omega_3) = \pm \frac{q {r}}{tk (q^2 + 1)^2}\big( \frac{2tk}{\sqrt{r}} \omega_1 \tensorc \omega_2 - \frac{2tk}{\sqrt{r}} \omega_2 \tensorc \omega_1 - \frac{t^2 k}{\sqrt{r}}\omega_3 \tensorc \omega_3 + tq \omega_3 \tensorc \omega_4 + tq \omega_4 \tensorc \omega_3 \big)
	\end{aligned}
\end{equation*}
Lastly, since $d(\omega_4) = 0$, $\nabla_0(\omega_4) = 0$\\
Thus, we are done with our proof. 
\end{proof}

\section{Existence of a unique bicovariant Levi-Civita connection} \label{section5}

In this section, we prove that except for finitely many $q \in (-1,1) \backslash \{ 0 \}$, the $4D_\pm$ calculi admit a unique bicovariant Levi-Civita connection for every bi-invariant pseudo-Riemannian metric (as defined in Definition \ref{24thmay20191}) on $\E$. We achieve this by verifying the hypotheses of Theorem \ref{15thjune20193}.

Recall that for the $4D_\pm$ calculus, we had the decomposition
\[ {}_0 \E \tensorc {}_0 \E = {\rm Ker}({}_0 \sigma - \id) \oplus {\rm Ker}({}_0 \sigma + q^2) \oplus {\rm Ker}({}_0 \sigma + q^{-2}). \]
We have already fixed the symbol ${}_0 \F$ for $ {\rm Ker}({}_0 \sigma + q^2) \oplus {\rm Ker}({}_0 \sigma + q^{-2}) $. Let us now denote ${\rm Ker}({}_0 \sigma - \id)$ by ${}_0 \E \tensorc^{\rm sym} {}_0 \E$. Moreover, as in Definition \ref{22ndaugust20192}, we define the $\IC$-linear map \[{}_0 (\Psym) : {}_0 \E \tensorc {}_0 \E \to {}_0 \E \tensorc {}_0 \E\]
to be the idempotent with range ${}_0 \E \tensorc^{\rm sym} {}_0 \E$ and kernel ${}_0 \F$. Since, ${}_0 (\Psym)$ is the idempotent onto the eigenspace of ${}_0 \sigma$ with eigenvalue one, and with kernel the eigenspaces with eigenvalues $q^2$ and $q^{-2}$, it is of the form (see (22) of \cite{article6})
\begin{equation} \label{27thjan20201}
{}_0(\Psym) = \frac{{}_0 \sigma + q^2}{1 + q^2}.\frac{{}_0 \sigma + q^{-2}}{1 + q^{-2}}.
\end{equation}. 

By Proposition \ref{26thjan20201}, the set $\{ \nu_i \}_{i=1}^{10}$ forms a basis of $\zeroE \tensorc^{\rm sym} \zeroE$, where $\nu_i$ are given as follows:
\begin{equation} \label{23rdjuly20191}
\begin{aligned} 
&\nu_1 = \omega_1 \tensorc \omega_1,
&\nu_2 =& \omega_2 \tensorc \omega_2,\\
&\nu_3 = \omega_3 \tensorc \omega_3 + t \omega_1 \tensorc \omega_2,
&\nu_4 =& \omega_4 \tensorc \omega_4,\\
&\nu_5 = \omega_2 \tensorc \omega_1 + \omega_1 \tensorc \omega_2,
&\nu_6 =& \omega_3 \tensorc \omega_2 + \frac{1}{q^2}\omega_2 \tensorc \omega_3,\\
&\nu_7 = \omega_3 \tensorc \omega_1 + q^2\omega_1 \tensorc \omega_3,
&\nu_8 =& \omega_4 \tensorc \omega_2 + \omega_2 \tensorc \omega_4 - \frac{t^2 k}{q^2 \sqrt{r}} \omega_2 \tensorc \omega_3,\\
&\nu_9 = \omega_4 \tensorc \omega_1 + \omega_1 \tensorc \omega_4 + \frac{t^2 k}{ \sqrt{r}} \omega_1 \tensorc \omega_3,
&\nu_{10} =& \omega_4 \tensorc \omega_3 + \omega_3 \tensorc \omega_4 + \frac{t^2 k}{q \sqrt{r}} \omega_1 \tensorc \omega_2.
\end{aligned}
\end{equation}
Thus, an arbitary element of $(\zeroE \tensorc^{\rm sym} \zeroE) \tensorc \zeroE$ is given by $X= \sum_{ij} A_{ij} \nu_i \tensorc \omega_j$, for some complex numbers $A_{ij}$. Hence, if we show that $({}_0(\Psym))_{23}(\sum_{ij} A_{ij} \nu_i \tensorc \omega_j )=0$ implies that $A_{ij}=0$ for all $i,j$, then $({}_0(\Psym))_{23}$ is a one-one map from $(\zeroE \tensorc^{\rm sym} \zeroE) \tensorc \zeroE$ to $\zeroE \tensorc (\zeroE \tensorc^{\rm sym} \zeroE)$. However, ${\rm dim}((\zeroE \tensorc^{\rm sym} \zeroE) \tensorc \zeroE) = {\rm dim}(\zeroE \tensorc (\zeroE \tensorc^{\rm sym} \zeroE))$, so that $({}_0(\Psym))_{23}$ is a vector space isomorphism from $(\zeroE \tensorc^{\rm sym} \zeroE) \tensorc \zeroE$ to $\zeroE \tensorc (\zeroE \tensorc^{\rm sym} \zeroE)$. Suppose $\{ A_{ij} \}_{ij}$ are complex numbers such that $({}_0(\Psym))_{23}(\sum_{ij} A_{ij} \nu_i \tensorc \omega_j) = 0$. Then, by \eqref{27thjan20201}, we have \begin{equation} \label{17thmarch2020} \big((q^2({}_0\sigma)_{23} + 1 )(({}_0\sigma)_{23} + q^2)\big)(\sum_{ij} A_{ij} \nu_i \tensorc \omega_j)=0. \end{equation}
We want to show that except for finitely many values of $q$, the above equation implies that all the $A_{ij}$ are equal to $0$. This involves a long computation, including a series of preparatory lemmas. We will be using the explicit form of ${}_0 \sigma(\omega_i \tensorc \omega_j)$ as given in Equation (4.1) of \cite{stachura} as well as \eqref{23rdjuly20191} to express the left hand side of \eqref{17thmarch2020} as a linear combination of basis elements $\omega_i \tensorc \omega_j \tensorc \omega_k$. Then we compare coefficients to derive relations among the $A_{ij}$. We do not provide the details of the computation. However, for the purposes of book-keeping, each equation is indexed by a triplet $(i,j,k)$ meaning that it is obtained by collecting coefficients of the basis element $\omega_i \tensorc \omega_j \tensorc \omega_k$ in the expansion of $\big((q^2({}_0\sigma)_{23} + 1 )(({}_0\sigma)_{23} + q^2)\big)(\sum_{mn} A_{mn} \nu_m \tensorc \omega_n)$.

\begin{lemma}
	We have the following equations:
	\begin{align} 
		\tag{1,1,1} &A_{11} =0\\
		\tag{1,1,2} &A_{12}(q^4 + 2) + (tA_{31} + A_{51} + A_{10,1} \frac{t^2 k}{q \sqrt{r}})2q^2 + (A_{73} q^2 + A_{93} \frac{t^2 k}{\sqrt{r}})2q(q^2-1) =0\\
		\tag{1,1,3} \begin{split}&A_{13}(q^4 + 2q^2 -1) + A_{14} (\frac{k}{\sqrt{r}}(q^2 -2 + q^{-2})) \\
			+ &(A_{71}q^2 + A_{91}\frac{t^2 k}{\sqrt{r}})2q^2 + A_{91}(\frac{k}{\sqrt{r}}q^{-2}(q^2 -1)) =0 \end{split}\\
		\tag{1,1,4} & A_{13}(-q^2 \frac{\sqrt{r}}{k}) + A_{14}(q^4 + 1) + (A_{71} q^2 + A_{91}\frac{t^2 k}{\sqrt{r}})\frac{\sqrt{r}}{k}q^4 + A_{91}(q^4 + 1)=0
	\end{align}
\end{lemma}
\begin{proof}
	The above equations are derived by comparing the coeffcients of $\omega_1 \tensorc \omega_1 \tensorc \omega_1$, $\omega_1 \tensorc \omega_1 \tensorc \omega_2$, $\omega_1 \tensorc \omega_1 \tensorc \omega_3$ and $\omega_1 \tensorc \omega_1 \tensorc \omega_4$ in \eqref{17thmarch2020}.
\end{proof}

\begin{lemma}
	We have the following equations:
	\begin{align}
		\tag{1,2,1} \begin{split} & A_{12}(2q^2 -1) + (t A_{31} + A_{51} + A_{10,1}\frac{t^2 k}{q \sqrt{r}})(q^4 + 1) + (A_{73} q^2 + A_{93} \frac{t^2 k}{\sqrt{r}})(-2 q (q^2 -1)) \\
			+ & A_{93}(- \frac{k}{q \sqrt{r}}(q^2 -1)^2) + (A_{74}q^2 + A_{94})(-\frac{k}{q\sqrt{r}}(q^2-1)^2)=0 \end{split}\\
		\tag{1,2,2} & t A_{32} + A_{52} + A_{10,2}\frac{t^2 k}{q \sqrt{r}}=0\\
		\tag{1,2,3} \begin{split} & (t A_{34} + A_{54} + A_{10,4}\frac{t^2 k}{q \sqrt{r}})(- \frac{k}{\sqrt{r}}(q^2 -1)^2) + (t A_{33} + A_{53} + A_{10,3}\frac{t^2 k}{q \sqrt{r}})(-(q^2-1)^2)\\
			+ & (A_{72} q^2 + A_{92}\frac{t^2 k}{\sqrt{r}})2q^2 + A_{92}(-\frac{k}{\sqrt{r}}(q^2 -1)^2)=0 \end{split}\\
		\tag{1,2,4} \begin{split} & (t A_{33} + A_{53} + A_{10,3}\frac{t^2 k}{q \sqrt{r}})\frac{q^4 \sqrt{r}}{k} + (t A_{34} + A_{54} + A_{10,4}\frac{t^2 k}{q \sqrt{r}})(q^4 +1)\\
			+ & (A_{72}q^2 + A_{92}\frac{t^2 k}{\sqrt{r}})(-q^2) + A_{92}(q^4+1) =0 \end{split}
	\end{align}
\end{lemma}
\begin{proof}
	The above equations are derived by comparing the coeffcients of $\omega_1 \tensorc \omega_2 \tensorc \omega_1$, $\omega_1 \tensorc \omega_2 \tensorc \omega_2$, $\omega_1 \tensorc \omega_2 \tensorc \omega_3$ and $\omega_1 \tensorc \omega_2 \tensorc \omega_4$ in \eqref{17thmarch2020}.
\end{proof}

\begin{lemma}
	We have the following equations:
	\begin{align}
		\tag{1,3,1} \begin{split} & A_{13}2 q^2 + A_{14} \frac{k}{\sqrt{r}}(-q^2 (q - q^{-1})^2) \\ + &(A_{71}q^2 + A_{91}\frac{t^2 k}{\sqrt{r}})(-q^4 + 2q^2 + 1) + A_{91}\frac{k}{\sqrt{r}}(-(q^2 -1)^2)=0 \end{split} \\
		\tag{1,3,2} \begin{split} & (t A_{33} + A_{53} + A_{10,3} \frac{t^2 k}{q \sqrt{r}})2 q^2 + (t A_{34} + A_{54} + A_{10,4} \frac{t^2 k}{q \sqrt{r}})\frac{k}{\sqrt{r}}(q^2 - 2q + q^{-2})\\ + & (A_{72} q^2 + A_{92} \frac{t^2 k}{\sqrt{r}})(q^4 + 2 q^2 -1) + A_{92}\frac{k}{\sqrt{r}}q^{-2}(q^2 - 1)^2 =0 \end{split} \\
		\tag{1,3,3} \begin{split} &(t A_{31} + A_{51} + A_{10,1}\frac{t^2 k}{q \sqrt{r}})(- 2 q^3 + 2 q) + A_{12}2 q (q^2 - 1) + A_{93}(- \frac{k}{\sqrt{r}}q^{-2}(q^2 - 1)^3 )\\ + & (A_{73} q^2 + A_{93} \frac{t^2 k}{\sqrt{r}})(- q^4 + 6 q^2 - 1) + (A_{74} q^2 + A_{94} \frac{t^2 k}{\sqrt{r}})(- \frac{k}{\sqrt{r}}q^{-2}(q^2 - 1)^3 ) = 0 \end{split} \\
		\tag{1,3,4} \begin{split} &(t A_{31} + A_{51} + A_{10,1}\frac{t^2 k}{q \sqrt{r}})(- 2 q^3 + 2 q) + (A_{73}q^2 + A_{93}\frac{t^2 k}{\sqrt{r}})\frac{\sqrt{r}}{k}q^4 \\ + & A_{93}(3(q^2 - 1)^2 + 2 q^2) + ( A_{74}q^2 + A_{94} \frac{t^2 k}{\sqrt{r}})(q^4 + 1) =0 \end{split}
	\end{align}
\end{lemma}
\begin{proof}
	The above equations are derived by comparing the coeffcients of $\omega_1 \tensorc \omega_3 \tensorc \omega_1$, $\omega_1 \tensorc \omega_3 \tensorc \omega_2$, $\omega_1 \tensorc \omega_3 \tensorc \omega_3$ and $\omega_1 \tensorc \omega_3 \tensorc \omega_4$ in \eqref{17thmarch2020}..
\end{proof}

\begin{lemma}
	We have the following equations:
	\begin{align}
		\tag{1,4,1} & A_{13}(-\frac{q^2 \sqrt{r}}{k}) + A_{14}(q^4 +1) + (A_{71}q^2 + A_{91}\frac{t^2 k}{\sqrt{r}})\frac{q^4 \sqrt{r}}{k} + A_{91}(q^4 +1)=0\\
		\tag{1,4,2} \begin{split} & (t A_{33} + A_{53} + A_{10,3}\frac{t^2 k}{q \sqrt{r}}) q^4 \frac{\sqrt{r}}{k} + (t A_{34} + A_{54} + A_{10,4}\frac{t^2 k}{q \sqrt{r}})(q^4 + 1) \\ + & ( A_{72}q^2 + A_{92} \frac{t^2 k}{\sqrt{r}})\frac{\sqrt{r}}{k}(- q^2) + A_{92}(q^4 + 1) =0 \end{split} \\
		\tag{1,4,3} \begin{split} & A_{12}(- \frac{r}{k} q^3) + (t A_{31} + A_{51} + A_{10,1}\frac{t^2 k}{q \sqrt{r}})\frac{\sqrt{r}}{k}q^3 + A_{93}(q^4 - 1)\\ + & (A_{73}q^2 + A_{93}\frac{t^2 k}{\sqrt{r}})\frac{\sqrt{r}}{k}q^2(q^2 - 1) + (A_{74} q^2 + A_{94} \frac{t^2 k}{\sqrt{r}})(q^4 + 1) = 0 \end{split} \\
		\tag{1,4,4} & A_{94}=0
	\end{align}
\end{lemma}
\begin{proof}
	The above equations are derived by comparing the coeffcients of $\omega_1 \tensorc \omega_4 \tensorc \omega_1$, $\omega_1 \tensorc \omega_4 \tensorc \omega_2$, $\omega_1 \tensorc \omega_4 \tensorc \omega_3$ and $\omega_1 \tensorc \omega_4 \tensorc \omega_4$ in \eqref{17thmarch2020}..
\end{proof}

\begin{lemma}
	We have the following equations:
	\begin{align} 
		\tag{2,1,1} &A_{51}=0\\
		\tag{2,1,2} \begin{split} & A_{52} (q^4 + 2) + A_{21} (2 q^2) + (A_{63} q^{-2} + A_{83} \frac{t^2 k}{q^2 \sqrt{r}})2q(q^2 -1) \\ + & A_{83}(\frac{k}{\sqrt{r}}q^{-1}(q^2 - 1)^2) + (A_{64}q^{-2} + A_{84} \frac{t^2 k}{q^2 \sqrt{r}})\frac{k}{\sqrt{r}}q(q^2 -2 + q^{-2}) = 0 \end{split} \\
		\tag{2,1,3} \begin{split} & A_{53}(q^4 + 2 q^2 -1) + A_{54}\frac{k}{\sqrt{r}}(q^2 -2 + q^{-2}) \\ + & (A_{64}q^{-2} + A_{84} \frac{t^2 k}{q^2 \sqrt{r}})2 q^2 + A_{81} \frac{k}{\sqrt{r}}q^{-2}(q^2 - 1) = 0 \end{split} \\
		\tag{2,1,4} & A_{53}(-\frac{\sqrt{r}}{k}q^2) + A_{54}(q^4 + 1) + (A_{61} q^{-2} + A_{81} \frac{t^2 k}{q^2 \sqrt{r}}) \frac{\sqrt{r}}{k}q^4 + A_{81}(q^4 + 1)= 0
	\end{align}
\end{lemma}
\begin{proof}
	The above equations are derived by comparing the coeffcients of $\omega_2 \tensorc \omega_1 \tensorc \omega_1$, $\omega_2 \tensorc \omega_1 \tensorc \omega_2$, $\omega_2 \tensorc \omega_1 \tensorc \omega_3$ and $\omega_2 \tensorc \omega_1 \tensorc \omega_4$ in \eqref{17thmarch2020}.
\end{proof}

\begin{lemma}
	We have the following equations:
	\begin{align}
		\tag{2,2,1} \begin{split} & A_{52}(2 q^2 - 1) + A_{21} (q^4 + 1) + (A_{63}q^{-2} + A_{83} \frac{t^2 k}{q^2 \sqrt{r}})(-2q(q^2 -1)) \\ + & A_{83}( \frac{k}{q \sqrt{r}}(q^2 -1)^2) + (A_{64}q^{-2} + A_{84} \frac{t^2 k}{q^2 \sqrt{r}})(-\frac{k}{\sqrt{r}}q(q^2 -2 + q^{-2})) =0 \end{split} \\
		\tag{2,2,2} &A_{22}=0\\
		\tag{2,2,3} \begin{split} & A_{23}(-q^4 + 2 q^2 + 1) + A_{24}(-\frac{k}{\sqrt{r}}(q^4 -2 q^2 + 1)) \\ + & (A_{62}q^{-2} + A_{82} \frac{t^2 k}{q^2 \sqrt{r}})2 q^2 + A_{82}(-\frac{k}{\sqrt{r}}(q^2 - 1)^2) =0 \end{split} \\
		\tag{2,2,4} & A_{23}\frac{q^4 \sqrt{r}}{k} + A_{24}(q^4 + 1) + (A_{62}q^{-2} + A_{82}\frac{t^2 k}{q^2 \sqrt{r}})\frac{\sqrt{r}}{k}(-q^2) + A_{82}(q^4 + 1) =0
	\end{align}
\end{lemma}
\begin{proof}
	The above equations are derived by comparing the coeffcients of $\omega_2 \tensorc \omega_2 \tensorc \omega_1$, $\omega_2 \tensorc \omega_2 \tensorc \omega_2$, $\omega_2 \tensorc \omega_2 \tensorc \omega_3$ and $\omega_2 \tensorc \omega_2 \tensorc \omega_4$ in \eqref{17thmarch2020}.
\end{proof}

\begin{lemma}
	We have the following equations:
	\begin{align}
		\tag{2,3,1} \begin{split} & A_{53}( 2 q^2) + A_{54}(q^4 + 1) + (A_{61}q^{-2} + A_{81}\frac{t^ k}{q^2 \sqrt{r}})(- q^4 + 2 q^2 + 1) \\ + & A_{81}\frac{k}{\sqrt{r}}(-(q^2 - 1)^2)=0 \end{split} \\
		\tag{2,3,2} \begin{split} & A_{23}2q^2 + A_{24}\frac{k}{\sqrt{r}}(q^2 - 2 + q^{-2}) + (A_{62}q^{-2} + A_{82}\frac{t^2 k}{q^2 \sqrt{r}})(q^4 + 2 q^2 -1) \\ + & A_{82}\frac{k}{\sqrt{r}}q^{-2}(q^2 - 1)^2 = 0 \end{split} \\
		\tag{2,3,3} \begin{split} & A_{52}2q(q^2 -1) + A_{21}(-2q^3 + 2q) + (A_{63 }q^{-2} + A_{83}\frac{t^2 k}{q^2 \sqrt{r}})(-q^4 + 6q^2 -1) \\ + & A_{83}\frac{k}{\sqrt{r}}(-q^{-2}(q^2 - 1)^3) + (A_{64}q^{-2} + A_{84}\frac{t^2 k}{q^2 \sqrt{r}})\frac{k}{\sqrt{r}}(-q(q - q^{-1})^3) = 0 \end{split} \\
		\tag{2,3,4} \begin{split} & A_{21}\frac{\sqrt{r}}{k}q^3 + (A_{63 }q^{-2} + A_{83}\frac{t^2 k}{q^2 \sqrt{r}})\frac{\sqrt{r}}{k}q^4 + A_{83}(3(q^2 - 1)^2 + 2 q^2) \\ + & (A_{64} q^{-2} + A_{84}\frac{t^2 k}{q^2 \sqrt{r}})(q^4 + 1) = 0 \end{split}
	\end{align}
\end{lemma}
\begin{proof}
	The above equations are derived by comparing the coeffcients of $\omega_2 \tensorc \omega_3 \tensorc \omega_1$, $\omega_2 \tensorc \omega_3 \tensorc \omega_2$, $\omega_2 \tensorc \omega_3 \tensorc \omega_3$ and $\omega_2 \tensorc \omega_3 \tensorc \omega_4$ in \eqref{17thmarch2020}.
\end{proof}

\begin{lemma}
	We have the following equations:
	\begin{align}
		\tag{2,4,1} & A_{53}\frac{\sqrt{r}}{k}(- q^2) + A_{54}(q^4 + 1) + (A_{61} q^{-2} + A_{81} \frac{t^2 k}{q^2 \sqrt{r}})\frac{\sqrt{r}}{k}q^4 + A_{81}(q^4 + 1) = 0\\
		\tag{2,4,2} & A_{23}\frac{\sqrt{r}}{k}q^4 + A_{24}(q^4 + 1) + (A_{62}q^{-2} + A_{82} \frac{t^2 k}{q^2 \sqrt{r}})\frac{\sqrt{r}}{k}( - q^2) + A_{82} (q^4 + 1) = 0\\
		\tag{2,4,3} \begin{split} & A_{52}\frac{\sqrt{r}}{k}(- q^3) + A_{21}\frac{\sqrt{r}}{k}q^3 + (A_{63 }q^{-2} + A_{83}\frac{t^2 k}{q^2 \sqrt{r}})\frac{\sqrt{r}}{k}q^2(q^2 -1) \\ + & A_{83}(q^4 - 1) + (A_{64}q^{-2} + A_{84}\frac{t^2 k}{q^2 \sqrt{r}})(q^4 + 1) = 0 \end{split} \\
		\tag{2,4,4} &A_{84}=0
	\end{align}
\end{lemma}
\begin{proof}
	The above equations are derived by comparing the coeffcients of $\omega_2 \tensorc \omega_4 \tensorc \omega_1$, $\omega_2 \tensorc \omega_4 \tensorc \omega_2$, $\omega_2 \tensorc \omega_4 \tensorc \omega_3$ and $\omega_2 \tensorc \omega_4 \tensorc \omega_4$ in \eqref{17thmarch2020}.
\end{proof}

\begin{lemma}
	We have the following equations:
	\begin{align} 
		\tag{3,1,1} &A_{71}=0\\
		\tag{3,1,2} \begin{split} & A_{72}(q^4 + 2) + A_{61}2q^2 + A_{33}2q(q^2 - 1) \\ + & A_{10,3}\frac{k}{\sqrt{r}}q^{-1}(q^2 - 1)^2 + A_{34} \frac{k}{\sqrt{r}}q(q^2 -2 + q^{-2})=0 \end{split} \\
		\tag{3,1,3} & A_{73}(q^4 + 2 q^2 - 1) + A_{74}\frac{k}{\sqrt{r}}(q^2 -2 + q^{-2}) + A_{31}2 q^2 + A_{10,1}\frac{k}{\sqrt{r}}q^{-2}(q^2 - 1) = 0\\
		\tag{3,1,4} & A_{73} \frac{\sqrt{r}}{k}(-q^{-2}) + A_{74} (q^4 + 1) + A_{31} \frac{\sqrt{r}}{k} q^4 + A_{10,1}(q^4 + 1) = 0
	\end{align}
\end{lemma}
\begin{proof}
	The above equations are derived by comparing the coeffcients of $\omega_3 \tensorc \omega_1 \tensorc \omega_1$, $\omega_3 \tensorc \omega_1 \tensorc \omega_2$, $\omega_3 \tensorc \omega_1 \tensorc \omega_3$ and $\omega_3 \tensorc \omega_1 \tensorc \omega_4$ in \eqref{17thmarch2020}.
\end{proof}

\begin{lemma}
	We have the following equations:
	\begin{align}
		\tag{3,2,1} \begin{split} & A_{72} (2 q^2 - 1) + A_{61}(q^4 + 1) + A_{33}2q (- (q^2 - 1)) \\ + & A_{10,33}\frac{k}{\sqrt{r}} (- q^{-1}(q^2 - 1)^2) + A_{34} \frac{k}{\sqrt{r}}(- q (q^2 - 2 + q^{-2})) = 0 \end{split} \\
		\tag{3,2,2} &A_{62}=0\\
		\tag{3,2,3} & A_{63}(-q^4 + 2 q^2 + 1) + A_{64}\frac{k}{\sqrt{r}}(-(q^4 - 2q^2 + 1)) + A_{32}2q^2 + A_{10,2}\frac{k}{\sqrt{r}}(-(q^2 - 1)^2) = 0\\
		\tag{3,2,4} & A_{63}\frac{\sqrt{r}}{k}q^4 + A_{64(q^4 + 1)} + A_{32}\frac{\sqrt{r}}{k}(- q^2) + A_{10,2}(q^4 + 1) = 0
	\end{align}
\end{lemma}
\begin{proof}
	The above equations are derived by comparing the coeffcients of $\omega_3 \tensorc \omega_2 \tensorc \omega_1$, $\omega_3 \tensorc \omega_2 \tensorc \omega_2$, $\omega_3 \tensorc \omega_2 \tensorc \omega_3$ and $\omega_3 \tensorc \omega_2 \tensorc \omega_4$ in \eqref{17thmarch2020}.
\end{proof}

\begin{lemma}
	We have the following equations:
	\begin{align}
		\tag{3,3,1} & A_{73}2q^2 + A_{74}\frac{k}{\sqrt{r}}(-(q^2 - 1)^2) + A_{31}(- q^4 + 2q^2 + 1) + A_{10,1} \frac{k}{\sqrt{r}}(-(q^2 - 1)^2) = 0\\
		\tag{3,3,2} & A_{63}2q^2 + A_{64}\frac{k}{\sqrt{r}}(q^2 - 2 + q^{-2}) + A_{32}(q^4 + 2q^2 - 1) + A_{10,2}\frac{k}{\sqrt{r}}q^{-2}(q^2 - 1)^2 = 0\\
		\tag{3,3,3} \begin{split} & A_{61} (-2q^3 + 2q) + A_{33} (- q^4 + 6 q^2 - 1) \\ + & A_{10,3}\frac{k}{\sqrt{r}}(- q^{-2} (q^2 - 1)^3) + A_{34} \frac{k}{\sqrt{r}}(-q (q - q^{-1})^3)=0 \end{split} \\
		\tag{3,3,4} & A_{61} \frac{\sqrt{r}}{k}q^3 + A_{33}\frac{\sqrt{r}}{k}q^4 + A_{10,3}(3(q^2 - 1)^2 + 2q^2) + A_{34}(q^4 + 1) = 0
	\end{align}
\end{lemma}
\begin{proof}
	The above equations are derived by comparing the coeffcients of $\omega_3 \tensorc \omega_3 \tensorc \omega_1$, $\omega_3 \tensorc \omega_3 \tensorc \omega_2$, $\omega_3 \tensorc \omega_3 \tensorc \omega_3$ and $\omega_3 \tensorc \omega_3 \tensorc \omega_4$ in \eqref{17thmarch2020}.
\end{proof}

\begin{lemma}
	We have the following equations:
	\begin{align}
		\tag{3,4,1} &A_{73}\frac{\sqrt{r}}{k}(-q^2) + A_{74}(q^4 + 1) + A_{31}\frac{\sqrt{r}}{k}q^4 + A_{10,1} (q^4 + 1) = 0\\
		\tag{3,4,2} & A_{63}\frac{\sqrt{r}}{k}q^4 + A_{64}(q^4 + 1) + A_{32}\frac{\sqrt{r}}{k}(-q^2) + A_{10,2}(q^4 + 1)=0\\
		\tag{3,4,3} & A_{72}\frac{\sqrt{r}}{k}(-q^3) + A_{61}\frac{\sqrt{r}}{k}q^3 = 0\\
		\tag{3,4,4} &A_{10,4}=0
	\end{align}
\end{lemma}
\begin{proof}
	The above equations are derived by comparing the coeffcients of $\omega_3 \tensorc \omega_4 \tensorc \omega_1$, $\omega_3 \tensorc \omega_4 \tensorc \omega_2$, $\omega_3 \tensorc \omega_4 \tensorc \omega_3$ and $\omega_3 \tensorc \omega_4 \tensorc \omega_4$ in \eqref{17thmarch2020}.
\end{proof}

\begin{lemma}
	We have the following equations:
	\begin{align} 
		\tag{4,1,1} &A_{91}=0\\
		\tag{4,1,2} \begin{split} & A_{92}(q^4 + 2) + A_{81}2q^2 + A_{10,3}2q(q^2 - 1) \\ + & A_{43} \frac{k}{\sqrt{r}}q^{-1}(q^2 - 1)^2 + A_{10,4}\frac{k}{\sqrt{r}}q(q^2 - 2 + q^{-2}) = 0 \end{split} \\
		\tag{4,1,3} & A_{93}(q^4 + 2q^2 - 1) + A_{94}\frac{k}{\sqrt{r}}(q^2 - 2 + q^{-2}) + A_{10,1}2q^2 + A_{41}\frac{k}{\sqrt{r}}q^{-2}(q^2 - 1) = 0\\
		\tag{4,1,4} & A_{93} \frac{\sqrt{r}}{k}(- q^2) + A_{94}(q^4 + 1) + A_{10,1} \frac{\sqrt{r}}{k}q^4 + A_{41}(q^4 + 1) = 0
	\end{align}
\end{lemma}
\begin{proof}
	The above equations are derived by comparing the coeffcients of $\omega_4 \tensorc \omega_1 \tensorc \omega_1$, $\omega_4 \tensorc \omega_1 \tensorc \omega_2$, $\omega_4 \tensorc \omega_1 \tensorc \omega_3$ and $\omega_4 \tensorc \omega_1 \tensorc \omega_4$ in \eqref{17thmarch2020}.
\end{proof}

\begin{lemma}
	We have the following equations:
	\begin{align}
		\tag{4,2,1} \begin{split} & A_{92}(2q^2 - 1) + A_{81}(q^4 + 1) + A_{10,3}2q(q^2 - 1) + A_{43} \frac{k}{\sqrt{r}} (- q^{-1} (q^2 - 1)^2) \\ + & A_{10,4} \frac{k}{\sqrt{r}}q(q^2 - 2 + q^{-2}) =0 \end{split} \\
		\tag{4,2,2} &A_{82}=0\\
		\tag{4,2,3} &A_{83} ( - q^4 + 2 q^2 + 1) + A_{84}\frac{k}{\sqrt{r}}(-q^4 + 2 q^2 - 1) + A_{10,2}2q^2 + A_{42}\frac{k}{\sqrt{r}}(q^2 - 1)^2 = 0\\
		\tag{4,2,4} & A_{83}\frac{\sqrt{r}}{k} q^4 + A_{84}(q^4 + 1) + A_{10,2}\frac{\sqrt{r}}{k}(- q^2) + A_{42}(q^4 + 1) = 0
	\end{align}
\end{lemma}
\begin{proof}
	The above equations are derived by comparing the coeffcients of $\omega_4 \tensorc \omega_2 \tensorc \omega_1$, $\omega_4 \tensorc \omega_2 \tensorc \omega_2$, $\omega_4 \tensorc \omega_2 \tensorc \omega_3$ and $\omega_4 \tensorc \omega_2 \tensorc \omega_4$ in \eqref{17thmarch2020}.
\end{proof}

\begin{lemma}
	We have the following equations:
	\begin{align}
		\tag{4,3,1} \begin{split} & A_{93} (q^4 + 2 q^2 - 1) + A_{94}\frac{k}{\sqrt{r}}(q^2 - 2 + q^{-2}) + A_{10,1}( - q^4 + 2q^2 + 1) \\ + & A_{41}\frac{k}{\sqrt{r}}(-(q^2 - 1)^2) = 0 \end{split} \\
		\tag{4,3,2} & A_{83}2q^2 + A_{84}\frac{k}{\sqrt{r}}(q^2 - 2 +q^{-2}) + A_{10,2}(q^4 + 2q^2 - 1) + A_{42}\frac{k}{\sqrt{r}}q^{-2}(q^2 - 1)^2=0\\
		\tag{4,3,3} \begin{split} & A_{92}2q(q^2 - 1) + A_{81}(-2q^3 + 2q) + A_{10,3}(-q^4 + 6q^2 - 1) + A_{43} \frac{k}{\sqrt{r}}(-q^{-2}(q^2 - 1)^3)\\ + &A_{10,4}\frac{k}{\sqrt{r}}(-q(q - q^{-1})^3) = 0 \end{split} \\
		\tag{4,3,4} & A_{81}\frac{\sqrt{r}}{k}q^3 + A_{10,3} \frac{\sqrt{r}}{k}q^4 + A_{43}( 3(q^2 - 1)^2 + 2q^2) + A_{10,4}(q^4 + 1) = 0
	\end{align}
\end{lemma}
\begin{proof}
	The above equations are derived by comparing the coeffcients of $\omega_4 \tensorc \omega_3 \tensorc \omega_1$, $\omega_4 \tensorc \omega_3 \tensorc \omega_2$, $\omega_4 \tensorc \omega_3 \tensorc \omega_3$ and $\omega_4 \tensorc \omega_3 \tensorc \omega_4$ in \eqref{17thmarch2020}.
\end{proof}

\begin{lemma}
	We have the following equations:
	\begin{align}
		\tag{4,4,1} & A_{93}\frac{\sqrt{r}}{k}(-q^2) + A_{94}(q^4 + 1) + A_{10,1}\frac{\sqrt{r}}{k}q^4 + A_{41}(q^4 + 1) = 0\\
		\tag{4,4,2} & A_{83}\frac{\sqrt{r}}{k}q^4 + A_{83}(q^4 + 1) + A_{10,2} \frac{\sqrt{r}}{k}(-q^2) + A_{42}(q^4 + 1) = 0\\
		\tag{4,4,3} &A_{92} \frac{\sqrt{r}}{k}(-q^3) + A_{81} \frac{\sqrt{r}}{k} q^3 + A_{10,3} \frac{\sqrt{r}}{k} q^2 (q^2 - 1) + A_{43} (q^4 - 1) + A_{10,4}(q^4 + 1)=0\\
		\tag{4,4,4} &A_{44}=0
	\end{align}
\end{lemma}
\begin{proof}
	The above equations are derived by comparing the coeffcients of $\omega_4 \tensorc \omega_4 \tensorc \omega_1$, $\omega_4 \tensorc \omega_4 \tensorc \omega_2$, $\omega_4 \tensorc \omega_4 \tensorc \omega_3$ and $\omega_4 \tensorc \omega_4 \tensorc \omega_4$ in \eqref{17thmarch2020}.
\end{proof}

\begin{thm} \label{14thfeb20201}
	For the $4D_\pm$ calculi, the map \[ ({}_0(\Psym))_{23}: (\zeroE \tensorc^{\rm sym} \zeroE) \tensorc \zeroE \rightarrow \zeroE \tensorc (\zeroE \tensorc^{\rm sym} \zeroE)  \] is an isomorphism except for, possibly, finitely many values of $q \in ( -1, 1 ) \backslash \{ 0 \}$. Hence, for each bi-invariant pseudo-Riemannian metric $g$, there exists a unique bicovariant Levi-Civita connection for each calculus.
\end{thm}
\begin{proof}
By the discussion preceding the above series of preparatory lemmas, we need to show that the system of equations given above admit only the trivial solution for $A_{ij}$, $i= 1, \dots, 10$, $j= 1, \dots, 4$. We then proceed to solve these equations for all $A_{ij}$. Note that the following variables are all identically zero in the above over-determined system:\\
$A_{11}$ (by (1,1,1)), $A_{94}$ (by (1,4,4)), $A_{51}$ (by (2,1,1)), $A_{22}$ (by (2,2,2)), $A_{84}$ (by (2,4,4)), $A_{71}$ (by (3,1,1)), $A_{62}$ (by (3,2,2)), $A_{10,4}$ (by (3,4,4)), $A_{91}$ (by (4,1,1)), $A_{82}$ (by (4,2,2)) and $A_{44}$ (by (4,4,4)).\\
This reduces the equations (1,3,1) and (1,4,1) to the following exact system of linear equations in the variables $A_{13}$ and $A_{14}$, with the associated matrix having determinant $q^2(q^2+1)^2$:
\begin{eqnarray*}
	 & A_{13}2 q^2 + A_{14} \frac{k}{\sqrt{r}}(-q^2 (q - q^{-1})^2) = 0\\
	 & A_{13}(-\frac{q^2 \sqrt{r}}{k}) + A_{14}(q^4 +1)  = 0
\end{eqnarray*}
Hence the solution for the variables $A_{13}$ and $A_{14}$ is zero.\\
We repeat this process for the rest of the $A_{ij}$, identifying a subset of equations which has been reduced to an exact one due to the previously solved $A_{ij}$, and then concluding that the set of $A_{ij}$ in the current set are also solved to be $0$ except for at most finitely many value of $q \in (-1,1) \backslash \{0\}$.\\
(2,2,3) and (2,2,4) reduce to the following system of linear equations in $A_{23}$ and $A_{24}$ with determinant $(q^2+1)^2$:\\
\begin{eqnarray*}
 & A_{23}(-q^4 + 2 q^2 + 1) + A_{24}(-\frac{k}{\sqrt{r}}(q^4 -2 q^2 + 1)) = 0 \\
 & A_{23}\frac{q^4 \sqrt{r}}{k} + A_{24}(q^4 + 1) = 0
\end{eqnarray*}
(4,1,3), (4,1,4) and (4,3,1) reduce to the following system of linear equations in $A_{41}$, $A_{93}$, $A_{10,1}$ with determinant $2q^{10} - 2q^4 - 2q^2 + 2$:\\
\begin{eqnarray*}
& A_{93}(q^4 + 2q^2 - 1) + A_{10,1}2q^2 + A_{41}\frac{k}{\sqrt{r}}q^{-2}(q^2 - 1) = 0 \\
& A_{93} \frac{\sqrt{r}}{k}(- q^2) + A_{10,1} \frac{\sqrt{r}}{k}q^4 + A_{41}(q^4 + 1) = 0 \\
& A_{93} (q^4 + 2 q^2 - 1) + A_{10,1}( - q^4 + 2q^2 + 1) + A_{41}\frac{k}{\sqrt{r}}(-(q^2 - 1)^2) = 0
\end{eqnarray*}
(4,1,2), (4,2,1), (4,3,3) and (4,4,3) reduce to the following system of linear equations in $A_{43}$, $A_{81}$, $A_{92}$, $A_{10,3}$ with determinant $4q^{14} + 10 q^{12} - 10 q^{10} - 8q^8 + 26 q^4 - 26 q^2 + 4$: \\
\begin{eqnarray*}
& A_{92}(q^4 + 2) + A_{81}2q^2 + A_{10,3}2q(q^2 - 1) + A_{43} \frac{k}{\sqrt{r}}q^{-1}(q^2 - 1)^2 = 0 \\
& A_{92}(2q^2 - 1) + A_{81}(q^4 + 1) + A_{10,3}2q(q^2 - 1) + A_{43} \frac{k}{\sqrt{r}} (- q^{-1} (q^2 - 1)^2) =0 \\
& A_{92}2q(q^2 - 1) + A_{81}(-2q^3 + 2q) + A_{10,3}(-q^4 + 6q^2 - 1) + A_{43} \frac{k}{\sqrt{r}}(-q^{-2}(q^2 - 1)^3) = 0 \\
& A_{92} \frac{\sqrt{r}}{k}(-q^3) + A_{81} \frac{\sqrt{r}}{k} q^3 + A_{10,3} \frac{\sqrt{r}}{k} q^2 (q^2 - 1) + A_{43} (q^4 - 1) = 0
\end{eqnarray*}
(3,4,3), (3,1,2), (3,2,1) and (3,3,3) reduce to the following system of linear equations in $A_{33}$, $A_{34}$, $A_{61}$, $A_{72}$ with determinant $-2q^2(q-1)^2(q+1)^2(q^2+1)^4$: \\
\begin{eqnarray*}
	&A_{72}\frac{\sqrt{r}}{k}(-q^3) + A_{61}\frac{\sqrt{r}}{k}q^3 = 0 \\
	&A_{72}(q^4 + 2) + A_{61}2q^2 + A_{33}(2q(q^2 - 1))  + A_{34} \frac{k}{\sqrt{r}}q(q^2 -2 + q^{-2})=0\\
	&A_{72} (2 q^2 - 1) + A_{61}(q^4 + 1) + A_{33} (- 2q(q^2 - 1)) + A_{34} \frac{k}{\sqrt{r}}(- q (q^2 - 2 + q^{-2})) = 0\\
	&A_{61} (-2q^3 + 2q) + A_{33} (- q^4 + 6 q^2 - 1) + A_{34} \frac{k}{\sqrt{r}}(-q (q - q^{-1})^3)=0
\end{eqnarray*}
(2,1,3) and (2,1,4) reduce to the following system of equations in $A_{53}$ and $A_{54}$ with determinant $q^4(q^2+1)^2$: \\
\begin{eqnarray*}
&A_{53}(q^4 + 2 q^2 -1) + A_{54}\frac{k}{\sqrt{r}}(q^2 -2 + q^{-2}) = 0\\
&A_{53}(-\frac{\sqrt{r}}{k}q^2) + A_{54}(q^4 + 1) = 0
\end{eqnarray*}
(1,1,2), (1,2,1), (1,3,3) and (1,3,4) reduce to a system of equations in $A_{12}$, $A_{31}$, $A_{73}$, $A_{74}$ with determinant a non-zero polynomial in $q$: \\
\begin{eqnarray*}
&A_{12}(q^4 + 2) + tA_{31} 2q^2 + A_{73} q^2 2q(q^2-1) =0\\
& A_{12}(2q^2 -1) + t A_{31} (q^4 + 1) + A_{73} q^2 (-2 q (q^2 -1)) \\
&t A_{31} (- 2 q^3 + 2 q) + A_{12}2 q (q^2 - 1) + A_{73} q^2 (- q^4 + 6 q^2 - 1) + A_{74} q^2 (- \frac{k}{\sqrt{r}}q^{-2}(q^2 - 1)^3 ) = 0 \\
&t A_{31} (- 2 q^3 + 2 q) + A_{73}q^2 \frac{\sqrt{r}}{k}q^4 +  A_{74}q^2(q^4 + 1) = 0
\end{eqnarray*}
(2,1,2), (2,2,1), (2,3,3), (2,3,4) and (2,4,3) reduce to a system of equations in $A_{21}$, $A_{52}$, $A_{63}$, $A_{64}$, $A_{83}$ with determinant a non-zero polynomial in $q$: \\
\begin{eqnarray*}
\begin{split} & A_{52} (q^4 + 2) + A_{21} (2 q^2) + (A_{63} q^{-2} + A_{83} \frac{t^2 k}{q^2 \sqrt{r}})2q(q^2 -1) \\ + & A_{83}(\frac{k}{\sqrt{r}}q^{-1}(q^2 - 1)^2) + (A_{64}q^{-2} + A_{84} \frac{t^2 k}{q^2 \sqrt{r}})\frac{k}{\sqrt{r}}q(q^2 -2 + q^{-2}) = 0 \end{split} \\
\begin{split} & A_{52}(2 q^2 - 1) + A_{21} (q^4 + 1) + (A_{63}q^{-2} + A_{83} \frac{t^2 k}{q^2 \sqrt{r}})(-2q(q^2 -1)) \\ + & A_{83}( \frac{k}{q \sqrt{r}}(q^2 -1)^2) + (A_{64}q^{-2} + A_{84} \frac{t^2 k}{q^2 \sqrt{r}})(-\frac{k}{\sqrt{r}}q(q^2 -2 + q^{-2})) =0 \end{split} \\
\begin{split} & A_{52}2q(q^2 -1) + A_{21}(-2q^3 + 2q) + (A_{63 }q^{-2} + A_{83}\frac{t^2 k}{q^2 \sqrt{r}})(-q^4 + 6q^2 -1) \\ + & A_{83}\frac{k}{\sqrt{r}}(-q^{-2}(q^2 - 1)^3) + (A_{64}q^{-2} + A_{84}\frac{t^2 k}{q^2 \sqrt{r}})\frac{k}{\sqrt{r}}(-q(q - q^{-1})^3) = 0 \end{split} \\
\begin{split} & A_{21}\frac{\sqrt{r}}{k}q^3 + (A_{63 }q^{-2} + A_{83}\frac{t^2 k}{q^2 \sqrt{r}})\frac{\sqrt{r}}{k}q^4 + A_{83}(3(q^2 - 1)^2 + 2 q^2) \\ + & (A_{64} q^{-2} + A_{84}\frac{t^2 k}{q^2 \sqrt{r}})(q^4 + 1) = 0
\end{split} \\
\begin{split} & A_{52}\frac{\sqrt{r}}{k}(- q^3) + A_{21}\frac{\sqrt{r}}{k}q^3 + (A_{63 }q^{-2} + A_{83}\frac{t^2 k}{q^2 \sqrt{r}})\frac{\sqrt{r}}{k}q^2(q^2 -1) \\ + & A_{83}(q^4 - 1) + (A_{64}q^{-2} + A_{84}\frac{t^2 k}{q^2 \sqrt{r}})(q^4 + 1) = 0 \end{split}
\end{eqnarray*}
(3,3,2) and (3,4,2) reduce to a system of equations in $A_{32}$, $A_{10,2}$ with determinant $q^4(q^2 + 1)^2$: \\
\begin{eqnarray*}
&  A_{32}(q^4 + 2q^2 - 1) + A_{10,2}\frac{k}{\sqrt{r}}q^{-2}(q^2 - 1)^2 = 0\\
&  A_{32}\frac{\sqrt{r}}{k}(-q^2) + A_{10,2}(q^4 + 1)=0\\
\end{eqnarray*}
Finally, (4,2,3) reduces identically to $A_{42} = 0$.\\
Hence we have shown that all $A_{ij}$ are identically equal to zero except for atmost finitely many values of $q \in (-1.1)$. Therefore, $({}_0(\Psym))_{23}|_{(\zeroE \tensorc^{\rm sym} \zeroE) \tensorc \zeroE}$ is an isomorphism if $q$ does not belong to this finite subset.\\
Since $SU_q(2)$ is a cosemisimple Hopf algebra, and we have shown that the map ${}_0 \sigma$ is diagonalisable, by Theorem \ref{15thjune20193}, for each bi-invariant pseudo-Riemannian metric $g$, each of the $4D_\pm$ calculi admits a unique bicovariant Levi-Civita connection for all but finitely many $q$.
\end{proof}

The proof of Theorem \ref{15thjune20193}, as given in \cite{article6}, involves explicitly constructing a Levi-Civita connection for each triple $(\E, d, g)$, subject to the accompanying hypothesis. In Theorem \ref{14thfeb20201}, we have shown that the hypothesis holds for the $4D_\pm$ calculi and for any bi-invariant pseudo-Riemannian metric. In this subsection, we provide the explicit construction of the Levi-Civita connection for a fixed arbitrary bi-invariant pseudo-Riemannian metric $g$. For this we will need to recall some definitions and results from \cite{article6}.

\begin{defn} %\label{24thaugust20195}
	Let $ \E$ and $ g $ be as above. We define a map 
	$$ V_g: {}_0\E \rightarrow ({}_0\E)^*, \qquad V_g (e) (f) = g (e \tensora f). $$
\end{defn}

\begin{defn} \label{27thaugust2019night1}
	Let $g$ be as above.
	We define a map 
	$$ g^{(2)}: ({}_0 \E \tensorc {}_0 \E) \tensorc ({}_0 \E \tensorc {}_0 \E) \rightarrow \mathbb{C} ~ {\rm by} ~ {\rm the} ~ {\rm formula} $$
	$$ g^{(2)} ((e_1 \tensorc e_2) \tensorc (e_3 \tensorc e_4)) = g(e_1 \tensora g(e_2 \tensora e_3) \tensora e_4) $$
	for all $e_1, e_2, e_3, e_4$ in ${}_0 \E.$
	
	We also define a map $ V_{g^{(2)}} : ({}_0 \E \tensorc {}_0 \E) \rightarrow ({}_0 \E \tensorc {}_0 \E)^*: = \Hom_{\mathbb{C}} ({}_0 \E \tensorc {}_0 \E, \mathbb{C}) $ by the formula
	$$ V_{g^{(2)}} (e_1 \tensorc e_2) (e_3 \tensorc e_4) = g^{(2)} ((e_1 \tensora e_2) \tensora (e_3 \tensora e_4)).$$
\end{defn}

\begin{prop}{(Propositions 4.4 and 4.9 of \cite{article6})} \label{14thfeb20202}
 The map $ V_{g} $ is one-one and hence a vector space isomorphism from $ {}_0 \E $ to $ ({}_0 \E)^*.$ Moreover, the map $V_{g^{(2)}}$ is a vector space isomorphism from ${}_0 \E \tensorc^{\rm sym} {}_0 \E$ onto $({}_0 \E \tensorc^{\rm sym} {}_0 \E)^*$.
\end{prop}

\begin{defn} \label{2ndaugust20191} 
	Let $V$ and $W$ be finite dimensional complex vector spaces. The canonical vector space isomorphism from $ V \tensora W^* $ to $ \Hom_\IC (W, V) $ will be denoted by the symbol $ \zeta_{V, W}. $ It is defined by the formula:
	\begin{equation} \label{24thaugust20194}
	\zeta_{V,W}(\sum_i v_i \tensora \phi_i)(w) = \sum_i v_i \phi_i(w).
	\end{equation}
\end{defn}

\begin{lemma}{(Lemma 3.12 of \cite{article6})} \label{15thjune20192}
	The following maps are vector space isomorphisms: 
	$$\zeta_{{}_0 \E \tensorc {}_0 \E , {}_0 \E}:({}_0 \E \tensorc^{\rm sym} {}_0 \E) \tensorc ({}_0 \E)^* \to {\rm Hom}_\IC({}_0 \E, {}_0 \E \tensorc^{\rm sym} {}_0 \E), $$ 
	$$ \zeta_{{}_0 \E, {}_0 \E \tensorc {}_0 \E }: {}_0 \E \tensorc ({}_0 \E \tensorc^{\rm sym} {}_0 \E)^* \to \homc ({}_0 \E \tensorc^{\rm sym} {}_0 \E, {}_0 \E).$$ 
\end{lemma}

\begin{defn} Given the maps $\zeta_{{}_0 \E \tensorc^{\rm sym} {}_0 \E , {}_0 \E}$, $\zeta_{{}_0 \E, {}_0 \E \tensorc^{\rm sym} {}_0 \E}$, $V_g$, $V_g^{(2)}$ and ${}_0(\Psym))_{23}$, the map \[\widetilde{\Phi_g}: \homc({}_0 \E, {}_0 \E \tensorc^{\rm sym} {}_0 \E) \to \homc({}_0 \E \tensorc^{\rm sym} {}_0 \E, {}_0 \E)\]
	is defined such that the following diagram commutes:
	\[
	\begin{tikzcd}
	\homc({}_0 \E, {}_0 \E \tensorc^{\rm sym} {}_0 \E) \arrow[r, rightarrow, "\zeta_{{}_0 \E \tensorc {}_0 \E , {}_0 \E}^{-1}"] \arrow[d, "\widetilde{\Phi_g}"] &[3 ex] ({}_0 \E \tensorc^{\rm sym} {}_0 \E) \tensorc ({}_0 \E)^* \arrow[r, "\id \tensorc V_g^{-1}"] & ({}_0 \E \tensorc^{\rm sym} {}_0 \E) \tensorc {}_0 \E \arrow[d, " ({}_0(\Psym))_{23} "]\\
	\homc({}_0 \E \tensorc^{\rm sym} {}_0 \E, {}_0 \E) &[3 ex] \arrow[ swap, l, " \zeta^{-1}_{{}_0 \E, {}_0 \E \tensorc {}_0 \E} "] {}_0 \E \tensorc ({}_0 \E \tensorc^{\rm sym} {}_0 \E)^* & \arrow[swap, l, "\id \tensorc V_{g^{(2)}}"] {}_0 \E \tensorc ({}_0 \E \tensorc^{\rm sym} {}_0 \E)
	\end{tikzcd} 
	\] 
\end{defn}

\begin{rmk}
	In Theorem \ref{14thfeb20201}, we proved that the map $ ({}_0(\Psym))_{23}: (\zeroE \tensorc^{\rm sym} \zeroE) \tensorc \zeroE \rightarrow \zeroE \tensorc (\zeroE \tensorc^{\rm sym} \zeroE)$	is an isomorphism. By Proposition \ref{14thfeb20202} and Lemma \ref{15thjune20192}, the remaining legs of the above commutative diagram are isomorphisms. Hence, $\widetilde{\Phi_g}: \homc({}_0 \E, {}_0 \E \tensorc^{\rm sym} {}_0 \E) \to \homc({}_0 \E \tensorc^{\rm sym} {}_0 \E, {}_0 \E) $ is also an isomorphism.
\end{rmk}

\begin{thm}
	For a fixed bi-invariant pseudo-Riemannian metric $g$, the bicovariant Levi-Civita connection $\nabla$ is defined on elements of ${}_0\E$ by
	\begin{equation} \label{14thfeb20203}
	\nabla = \nabla_0 + \widetilde{\Phi_g}^{-1}(dg - \widetilde{\Pi^0_g}(\nabla_0)),
	\end{equation}
	where $\nabla_0$ is the bicovariant torsionless connection constructed in Theorem \ref{31stmar20201}. Here, $\nabla_0$ and $g$ are considered as restrictions on ${}_0 \E$ and ${}_0 \E \tensorc {}_0 \E$ respectively.
\end{thm}
\begin{proof}
	Let us recall from Remark \ref{26thjan20202}, that it is sufficient to define a connection on ${}_0 \E$ to define it on the whole of $\E$. Next, by \eqref{1stmay20191}, $\widetilde{\Pi^0_g}(\nabla_0)$ is a well-defined map in $\homc({}_0\E \tensorc {}_0 \E, {}_0\E)$. The map $dg$ is a well-defined map in $\homc({}_0\E \tensorc {}_0 \E, {}_0\E)$. (Indeed it is the zero-map, since $g$ maps ${}_0\E \tensorc {}_0 \E$ to $\IC$, and $d$ maps $\IC$ to $0$. That we write it at all in the formula of $\nabla$ is because of how it appears in the proof of Proposition 7.3 of \cite{article6}.) We have already remarked that $\widetilde{\Phi_g}$ is a well-defined isomoprhism from $ \homc({}_0 \E, {}_0 \E \tensorc^{\rm sym} {}_0 \E)$ to $\homc({}_0 \E \tensorc {}_0 \E, {}_0 \E)$. Hence, the right-hand side of \eqref{14thfeb20203} is a well-defined map in $\homc({}_0 \E, {}_0 \E \tensorc {}_0 \E)$. That it defines the unique bicovariant Levi-Civita connection on $\E$ follows from the proofs of Proposition 7.3 and Theorem 7.8 of \cite{article6}, and we leave out the details.
\end{proof}

\end{document}